\documentclass[letterpaper,11pt,draft,runningheads]{llncs}
\usepackage{amssymb,theorem}
\usepackage{amsmath,array,amscd}
\usepackage[mathscr]{eucal}
\def\N{{\rm I}\hskip -.22em {\rm N}}

\def\F{{\rm I}\hskip -.22em {\rm F}}

\newtheorem{algo}{Algorithm}
\newcommand{{\M}}{{\mathcal N}}
\newcommand{{\Pic}}{{\rm Pic}}
\newcommand{{\Cl}}{{\rm Cl}}
\newcommand{\IN}{{\tt INPUT:\ \ }}
\newcommand{\OUT}{{\tt OUTPUT:\ }}

\title{On Using Expansions to the Base of $-2$}

\author{Roberto Avanzi\inst{1}, Gerhard Frey\inst{1}, Tanja Lange\inst{2}, Roger Oyono\inst{1}}

\institute{IEM, University of Duisburg-Essen\\
  Ellernstrasse 29, D-45326 Essen, Germany\\
 {\tt \{mocenigo,frey,oyono\}@exp-math.uni-essen.de}\\ \and
 ITSC, Ruhr-University of Bochum,\\ Universit\"atsstr. 150, D-44780 Bochum, Germany,\\{\tt
lange@itsc.ruhr-uni-bochum.de}}
\date{123o23}
\begin{document}


\maketitle
\begin{abstract}
This short note investigates the effects of using expansions to the
base of $-2$. The main applications we have in mind are cryptographic
protocols, where the crucial operation is computation of scalar
multiples.  For the recently proposed groups arising from Picard curves
this leads to a saving of
at least
7\% for the computation of an $m$-fold. For more general
non-hyperelliptic genus 3 curves we expect a larger speed-up.\\[.5ex]
{\em ACM Computing Reviews Categories:} E3 public key cryptosystems, G4 efficiency
\end{abstract}
\section{Introduction}
Recently, groups associated to 
elliptic and
hyperelliptic curves received a lot of
attention for cryptographic applications, and further kinds of curves
were proposed and their arithmetic studied intensively.
They allow smaller operands compared to RSA and DL in finite fields
making them attractive for restricted devices. The performance on
such small units is good~\cite{PeWoGuPa03}. More general
curves were suggested and the group operations optimized for
cryptographic applications.

To compute scalar multiples $mD$, binary expansions of $m$ are used
and the computation of $mD$ is split up as a sequence of additions and
doublings.  To achieve faster computations one uses windowing methods
and signed representations (for a broad overview see
Knuth~\cite{kn298}).

Our idea speeds up scalar multiplication in groups for which computing
$-2D$ and $-(D_1+D_2)$ is faster than computing $2D$ and $D_1+D_2$,
respectively.  For elliptic and hyperelliptic curves the negative of
an element can be obtained almost for free.  Hence, these groups will
most probably not benefit from our new idea.  But, the situation is
different for Picard curves or more general genus 3 quartic curves:
using the $-2$-adic expansion instead of a $2$-adic we reduce the
complexity by at least 7\%.

Of course, our considerations are not restricted to cryptography but
allow speeding up scalar multiplication in groups in which ``adding up
to the neutral element'' is easier computed than addition and so they
could be of use in computer algebra systems, too.

In this note we first describe the idea of $-2$-adic expansions and
show how to apply them.  Then we sketch the applications we have in
mind, and finally show the time saving for Picard curves.

To fix notation, let $G$ be a finite abelian group of order $\ell$ and
let $D$ be a generator of $G$. Furthermore, we assume that computing
$-2D$ or $-(D_1+D_2)$ is actually faster than computing $2D$ or
$D_1+D_2$ respectively.

\section{$-2$-adic Expansions}
Assume that we want to compute $mD,\ m< \ell$ and put $l(m)$ the
length and $w(m)$ the number of nonzero bits of the used expansion of
$m$.  Since we need both $D$ and $-D$, we can allow signed digit
representations.

If $l(m)+w(m)$ is even we start with $D$, otherwise with $-D$.  While
the doubling is always replaced by the computation of $-2$-times the
intermediate result $E$ we need to pay a little more attention on how
to perform the former additions and subtractions.

\begin{algo}
\label{algo}
\ \\
\IN $m\in \N$, $m=\sum_{i=0}^{l(m)-1}m_i2^i,\, m_i\in\{0,\pm 1\}$, $D\in G$\\
\OUT $E:=mD$
\begin{enumerate}
\item precompute and store $-D$
\item compute $l(m), w(m)$;
\item put $E:=(-1)^{f}D$,\ \ where\ \ $f:= l(m)+w(m)\bmod 2$;
\item\label{loop} for $i=l(m)-2$ to $0$ do
\begin{enumerate}
\item $E:=-2(E)$;
\item $f:=1 - f$;
\item if $m_i\neq 0$
\begin{enumerate}
\item \label{keineSorge}
$E:=-(E+(-1)^{f}m_iD)$;
\item $f:=1-f$;
\end{enumerate}
\end{enumerate}
\item output($E$);
\end{enumerate}
\end{algo}

The correctness follows from the fact that $l(m)+w(m)-2$ is the total
number of minus signs in front of the initial $E$. During the process
$f$ keeps track of the parity of the number of sign changes.  In step
{\it \ref{keineSorge}}, no multiplication is required. $f$ assumes only
values in $\{0,1\}$ and $m_i\in\{\pm 1\}$.

Using this idea introduces only little
bookkeeping overhead, namely the additional variable $f$.  So such a
system is really practicable -- and  useful if the operations
involving the negative signs are faster.

\medskip\penalty-300
{\bf Remarks:}\vspace*{-7.534281pt}%
\begin{enumerate}
\item If the expansion of $m$ is not computed beforehand, one can
always start with $E=D, f=0$, irrespective of the parity of
$l(m)+w(m)$.  The loop {\it \ref{loop}.} is performed as above.
Before Step {\it 5.}  one checks whether $f=1$, in which case one
outputs $-E$ instead.  So one avoids precomputing the expansion, at
the price of a second negation with probability 1/2.
\item If only for one of addition or doubling the negative is faster,
similar considerations hold if one only replaces that operation.
\item Of course the method can be combined with signed sliding
windowing methods.  In the applications we have in mind $rD$ and $-rD$
can be computed with only a few more operations than $rD$ 
alone.
\end{enumerate}

\section{Applications}
In this section we need to state some details from mathematics to show
that there actually are applications of our idea. For an introduction
to hyperelliptic curves see \cite{MWZohneK}. The following holds for
arbitrary curves.

Let $C$ be a curve of genus $g$. The group used for cryptographic
applications is a subgroup of the {\em divisor class group} of $C$:
we briefly recall its main properties.
Let $P_\infty \in C$ be fixed.
A divisor  is a formal sum of points.
We are interested in 
the degree zero divisors given by sums
\begin{equation}
D-nP_\infty=\sum_{i=0}^nP_i-nP_\infty\enspace,\quad P_i\in C\smallsetminus \{P_\infty\}\enspace.
\label{gaehn}
\end{equation}
The principal divisors are the divisors of functions.  The divisor
class group is the group of the degree zero divisors modulo the
principal ones.  In each divisor class there exists a unique element
\eqref{gaehn} with $n\leq g$ 
minimal.

To add two classes $c_1,\,c_2$ one formally adds the representing
divisors: $D_1+D_2-(n_1+n_2)P_\infty$.  Then one determines a function
$f$ passing through the points on $D_1+D_2$ with poles only in $k
P_\infty$ ($k$ minimal) with multiplicities taken into account. Let
$D_3$ be the divisor represented by the points of intersection of $f$
with the curve
which are 
not in $D_1$ and $D_2$. Put $n_3=k-n_1-n_2$ and let $c_3$ be the class
of $D_3-n_3P_\infty$.  Since $D_1+D_2+D_3 -k P_\infty$ add up to a
principal divisor we get $c_1+c_2+c_3=0$.

Usually one proceeds to find the negation $c_4=-c_3$ to get
$c_1+c_2=c_4$. Our new proposal allows to skip this last step.
Doublings just work the same with the function passing through the
points of the representing divisor with doubled multiplicity.

For hyperelliptic curves, taking the negative is very simple. 
The formulae for genus $2$ and $3$ \cite{TLexpljoint,pelzl}
reveal that computing $-2D$ instead of $2D$ saves only some additions
in the underlying field. Therefore, we do not expect the $-2$-adic
expansion to lead to a saving.

{\em But the situation is completely different for non-hyperelliptic
curves}. For char$(\F_q)\ne 3$ a Picard curve 
%
%
can be given by:
\[ z y^3 = z^4 f_4 (x/z)\enspace, \quad f_4\in \F_q[x]\]
where $f_4$ is  monic,  square-free and of degree 4.

The arithmetic on Picard curves is detailed in \cite{FO} (see also
\cite{BEFG}). An addition needs 144M, 12S, and 2I and a doubling 158M,
16S, and 2I in the generic case.  Applying Algorithm~\ref{algo}
reduces the costs to 133M, 9S, 2I or 
147M, 
13S, 2I respectively. Some
field additions are saved as well.  Thus, here the saving is at least
7.5\% or 7\%, respectively, assuming a ratio of $10:1$ for inversions
and $2:3$ for squarings in relation to multiplications.

An ordinary genus 3 quartic over $\F_q$ is given by a projective equation:
\[ \left(a_1x^2+a_2y^2+a_3z^2+a_4xy+a_5xz+a_6yz\right)^2=l(x,y,z)xyz,\]
where $l\in  \F_q[x,y,z]$ is linear and $a_i\in \F_q$.
On these curves computing the negation is even more
complicated than on Picard curves (see
\cite{FOR}). Therefore, the saving due to $-2$-adic
expansions is more dramatic.

\bibliographystyle{plain} \bibliography{../../bibfiles/bib}
\end{document}